\documentclass[11pt]{amsart}

\usepackage{bbm}
\usepackage[margin=1in]{geometry}
\usepackage{amsmath,amssymb,amsthm,mathtools}
\usepackage{enumitem}
\usepackage[colorlinks=true,linkcolor=blue,citecolor=blue,urlcolor=blue]{hyperref}
\usepackage[nameinlink,capitalize]{cleveref}
\usepackage{tikz-cd}

\newtheorem{theorem}{Theorem}[section]

\newtheorem{conjecture}[theorem]{Conjecture}
\newtheorem{question}[theorem]{Question}
\newtheorem{proposition}[theorem]{Proposition}
\newtheorem{lemma}[theorem]{Lemma}
\newtheorem{corollary}[theorem]{Corollary}

\theoremstyle{definition}

\newcommand{\Fq}{\mathbb F_q}

\newcommand{\N}{\mathbb{N}}

\newcommand{\PP}{\mathbb{P}}

\newcommand{\Wb}{\overline{W}}

\newcommand{\sg}{S(\Gamma)}
\newcommand{\ag}{A(\Gamma)}

\newcommand{\Supp}{\operatorname{Supp}}
\newcommand{\Fb}{\Bar{F}}

\newcommand{\itr}{{\langle i \rangle}}

\newcommand{\dimk}{\dim_{\kappa}}
\newcommand{\ka}{\kappa}

\newcommand{\Ann}{{\rm Ann}}

\newcommand{\Pm}{\mathbb{P}^m}

\newcommand{\gG}{g_{\Gamma}}
\newcommand{\Ag}{A(\Gamma)}

\newcommand{\ga}{g_{A(\Gamma)}}
\newcommand{\ggo}{g_{\Gamma_1}}
\newcommand{\go}{\Gamma_1}
\newcommand{\sgo}{S(\Gamma_1)}
\newcommand{\tsum}{\mathop{\textstyle\sum}\nolimits}

\newcommand{\Res}{\text{Res}}

\newcommand\restr[2]{{% we make the whole thing an ordinary symbol
  \left.\kern-\nulldelimiterspace % automatically resize the bar with \right
  #1 % the function
  \vphantom{\big|} % pretend it's a little taller at normal size
  \right|_{#2} % this is the delimiter
  }}

\DeclarePairedDelimiter\ceil{\lceil}{\rceil}

\title{Artinian reduction applied to the Largest Intersection Problem}
\author{Yuxin Lin}

\date{\today}

\begin{document}

\begin{abstract}
Suppose we have $r$ hypersurfaces in $\PP^m$ of degree $d$, whose defining polynomials are linearly independent and whose intersection is zero-dimensional. What is the maximum number of points in their intersection? This problem is motivated by computing the generalized Hamming weights of projective Reed–Muller codes, and a conjectured formula for this maximum was given in \cite{lin2026largest}. We prove the conjectured formula for arbitrary $m$ when $\tbinom{m+d}{d}-\binom{d+2}{2}<r\leq\tbinom{m+d}{d}$. We use the notion of Artinian reduction to reformulate the problem in terms of graded Artinian algebras and then apply combinatorial tools to bound their Hilbert functions.

% Our proof uses Artinian reduction to bound the failure of zero-dimensional schemes contained in a complete intersection to impose independent conditions on degree-$d$ polynomials.
\end{abstract}

\maketitle

\section{Introduction}

Let $\ka$ be an algebraically closed field (possibly of characteristic zero).
Let $S(m)=\ka[x_0,\dots,x_m]$ and $S(m)_d$ be its $d$-th graded component. 
Given a subset $W\subseteq S(m)_d$, we denote the common vanishing set by $V(W) \subseteq \PP^{m}$.
The following problem was raised in \cite{lin2026largest}.
\begin{question}\label{Question: zero dim problem}
Suppose we are given $m,d\geq 1$ and $m\leq r\leq \binom{m+d}{d}$. What is the largest $|V(W)|$ as $W$ varies through $r$-dimensional vector subspaces of $S(m)_d$ for which $\dim(V(W))=0$?
\end{question}
They denote this quantity as $u_r(d,m)$ and conjecture an exact formula for it.
We will introduce some notation to state their conjecture.
Let
\[
\Omega'(d,m)
:=\Big\{
(\alpha_1,\alpha_{2},\dots,\alpha_{m+1})\in \N_0^{m+1}
\mid \tsum_{i=1}^{m+1}\alpha_i=d, 
d\notin\{\alpha_2,\dots,\alpha_{m}\}  
\Big\}.
\]
So $|\Omega'(d,m)|=\binom{m+d}{d}-(m-1)$. Order $\Omega'(d,m)$ by lexicographical ordering and let $\omega'_r(d,m)$ be its $r$-th largest element. If $\omega'_r(d,m)=(\alpha_1,\dots,\alpha_{m+1})$, they define
\[
H'_{r}(d,m)
:=\alpha_1 d^{m-1}+\dots+\alpha_m.
\]

\begin{conjecture}{\cite[Conjecture 1.6]{lin2026largest}}\label{Conjecture urdm}
Given $d,m\geq 1$ and $m\leq r\leq \binom{d+m}{m}$, we have
\[
u_r(d,m)=H_{r-(m-1)}'(d,m).
\]
\end{conjecture}
They prove their conjecture in the following cases, notably for $m=2$.
\begin{theorem}{\cite[Theorem 1.8]{lin2026largest}}\label{thm: cases of urdm}
Given positive integers $d,m,r$ with $m\leq r\leq \binom{d+m}{m}$, if at least one of the following holds
\begin{enumerate}
    \item $d=1$;
    \item $m=1$ or $m=2$;
    % \item $m=2$;
    \item $r=m$;
    \item $\binom{m+d}{d}-d\leq r\leq \binom{m+d}{d}$,
\end{enumerate}
then we have $u_r(d,m)=H_{r-(m-1)}'(d,m)$.
\end{theorem}
In general, they prove that their conjectured formula is at least a lower bound.
\begin{proposition}{\cite[Proposition 1.7]{lin2026largest}}\label{Prop: lower bound urdm}
For positive integers $m,d,r$ with $m\leq r\leq \binom{d+m}{m}$, we have
\[
u_r(d,m)\geq H_{r-(m-1)}'(d,m).
\]
\end{proposition}

Our main result is to extend the conclusion of Theorem \ref{thm: cases of urdm} to $\binom{m+d}{d}-\tbinom{d+2}{2}< r\leq \tbinom{m+d}{d}$. 

\begin{theorem}\label{Thm: main result}
Given positive integers $d,m,r$ with
$m\leq r\leq \binom{m+d}{d}$ and
$r>\binom{m+d}{d}-\binom{d+2}{2}$, we have
\[
u_r(d,m)=H_{r-(m-1)}'(d,m).
\]
\end{theorem}
We prove Theorem~\ref{Thm: main result} by bounding the number of independent conditions imposed by zero-dimensional schemes on degree $d$ polynomials.

\subsection{Number of independent conditions}

Let $\Gamma$ be a zero-dimensional subscheme of $\PP^m$.
Let $I(\Gamma)$ be its homogeneous ideal, and let $S(\Gamma)=\kappa[x_0,\dots,x_m]/I(\Gamma)$ be its graded coordinate ring.
We denote their $k$-th homogeneous pieces as $I(\Gamma)_k$ and $S(\Gamma)_k$ respectively.
Note that $I(\Gamma)_k$ is the vector space of degree $k$ polynomials that vanish on $\Gamma$. 
The Hilbert function $h_{\sg}(k):=\dimk S(\Gamma)_k$ can also be viewed as the number of independent conditions imposed by $\Gamma$ on degree $k$ polynomials.
We denote by $g_{\Gamma}(k)$ the failure of $\Gamma$ to impose independent conditions on polynomials of degree $k$. That is,
\[
g_{\Gamma}(k)
:=\deg(\Gamma)-h_{\sg}(k) 
=\deg(\Gamma)+\dimk I(\Gamma)_k-\tbinom{m+k}{k}.
\]
Note that $g_{\Gamma}(d)$ can also be interpreted in terms of the cohomology of the ideal sheaf, $g_{\Gamma}(d)=h^1(\PP^m,\mathcal{I}_\Gamma(d))$.
Here, $\deg(\Gamma)$ is the degree of the zero-dimensional subscheme $\Gamma$ in $\PP^m$ as defined in \cite[Section 1.7]{hartshorne2013algebraic}. Since $\Gamma$ is zero-dimensional, it is the number of points in $\Gamma$ counted with multiplicity, denoted by $|\Gamma|_{\mathfrak{m}}$ in \cite[Definition 2.1]{lin2026largest}.
The problem of when a set of points imposes independent conditions, that is, when $g_{\Gamma}(d)$ is zero, has been widely studied \cite{ballico2014zero, ballico2013gaps, caviglia2021cayley, couvreur2012dual, gold2005cayley, lin2026largest, moreno2026generalized}.

% If $\Gamma$ is reduced, then $\deg(\Gamma)$ is just the number of points in $\Gamma$.
% \textcolor{blue}{Unify the notation: $\gG(k)$ failure to impose independent condition; $h_{S(\Gamma)}(k)$ for Hilbert function, $A(\Gamma)$ for the Artinian reduction, $I(\Gamma)_k$, $S(\Gamma)_k$ for the $k^{th}$ graded piece }
In \cite{lin2026largest}, the authors proved the case $m=2$ of Theorem~\ref{thm: cases of urdm} using the following bound.

\begin{proposition}{\cite[Lemma 5.2, Proposition 5.6]{lin2026largest}}\label{Prop: cond on gx udr}
Consider positive integers $t$, $d$ satisfying $2\leq t\leq d$.
Given coprime polynomials $G_1$, $G_2$ in $S(2)_d$ and a subscheme $\Gamma\subseteq V(G_1,G_2)$, with $\deg(\Gamma)\leq td+d-t+1$, we have $g_{\Gamma}(d)\leq \tbinom{t-1}{2}$.
\end{proposition}

In this paper, we will prove the following extension of Proposition \ref{Prop: cond on gx udr}.

\begin{theorem}\label{Thm: main result g_Gamma}
Let $G_1, \dots, G_m\in S(m)_d$ be a regular sequence. For every sub-scheme $\Gamma\subseteq V(G_1,\dots, G_m)$:
\begin{enumerate}
\item if $\deg(\Gamma) \leq 3d-1$, then $\gG(d)=0$;
\item for $2\leq t\leq d$, if $\deg(\Gamma)\leq td+d-t+1$, then $g_{\Gamma}(d)\leq \tbinom{t-1}{2}$;
\item if $\deg(\Gamma) \leq d^2+2$, then $\gG(d) \leq \binom{d-1}{2}$.
\end{enumerate}
\end{theorem}

Note that (3) is close to (2) with $t=d$, but it allows for an extra value of $\deg(\Gamma)=d^2+2$.
For $d\geq 2$, part (1) is the case $t=2$ of part (2).
We state it separately because it admits a generalization
that does not require $\Gamma$ to lie in a complete intersection.

\begin{theorem}\label{Thm: 3d-1 geometric}
Let $\Gamma\subseteq \PP^m$ be a zero-dimensional subscheme with $\deg(\Gamma)\leq 3d-1$. 
We have $g_{\Gamma}(d)=0$ if and only if
\begin{enumerate}
\item for every line $L$, $\deg(L\cap \Gamma)\leq d+1$;
\item for every planar conic $C$, $\deg(C\cap \Gamma)\leq 2d+1$.
\end{enumerate}
\end{theorem}
Couvreur \cite{couvreur2012dual} proved Theorem~\ref{Thm: 3d-1 geometric} for reduced schemes. Ballico \cite{ballico2013gaps} proved it in characteristic zero, with the additional requirement that $\deg(\Gamma)-\deg(\Gamma_{\textrm{red}})\leq d$ when $m\geq 3$.
Ellia and Peskine in \cite{ellia2006groupes} proved a result that implies Theorem~\ref{Thm: 3d-1 geometric} for $m=2$.

We will prove these results by connecting them to standard graded Artinian algebras. 

% \textcolor{blue}{I think we should include all the artinian formulation to section \ref{Subsec: Artinian reduction}, including this theorem. The main result is on zero-dim scheme.}

% \textcolor{red}{The following part goes away: We explain the notion of the Artinian reduction of a zero dimensional scheme $\Gamma$ in Section~\ref{Subsec: Artinian reduction}. This was first introduced by Eisenbud, Green and Harris in \cite{eisenbud1996cayley}.
% The following result is the Artinian algebra analog of Theorem~\ref{Thm: main result g_Gamma}.
% \begin{theorem}\label{Thm: main result algebra}
% Consider $2\leq t\leq d$.
% Suppose $J\subseteq \ka[x_1, \dots, x_m]$ is a graded ideal that contains an regular sequence of degree $(\underbrace{d,\ldots,d}_{m})$. 
% Let $A=\ka[x_1, \dots, x_m]/J$.
% If $\dim(A)\leq td+d-t+1$, then
% \[
% \dim(A_{\geq d+1})\leq \tbinom{t-1}{2}.
% \]
% \end{theorem}}

We give a brief outline of the paper.
In Section~\ref{Section: Background}, we explain how Question~\ref{Question: zero dim problem} relates to other problems studied in the literature. We also explain the notion of Artinian reduction, give an Artinian formulation of Theorem \ref{Thm: main result g_Gamma}, and give some lemmas that we will use later.
In Section~\ref{sec: reduction from 1.8 to 3.1}, we state a Artinian formulation of Theorem~\ref{Thm: 3d-1 geometric}, and show that the Artinian formulation implies Theorem~\ref{Thm: 3d-1 geometric}.
%Theorem~\ref{thm: 3d-1 artinian} implies Theorem~\ref{Thm: 3d-1 geometric}. 
In Section~\ref{sec: smallest and largest r}, we prove the Artinian formulations of Theorem~\ref{Thm: 3d-1 geometric} and parts (1) and (3) of Theorem~\ref{Thm: main result g_Gamma}. Finally, in Section~\ref{Section: Artinian tail} we prove the Artinian formulation of part (2) of Theorem~\ref{Thm: main result g_Gamma}.

% Proposition \ref{Prop: algebra reg seq 3d-1} and \ref{prop: end bound}.

\section{Background}\label{Section: Background}

\subsection{Beelen-Datta-Ghorpade conjecture and $k$-dimensional conjecture}

Question \ref{Question: zero dim problem} was motivated by the following conjecture of Beelen, Datta and Ghorpade. 

Let $m,d\ge1$, $q\ge d+1$ be a prime power, and $1\le r\le\binom{m+d}{d}$. For an \(r\)-dimensional vector subspace \(W\subseteq\mathbb F_q[x_0,\ldots,x_m]_d\), denote by \(V(W)\subseteq\mathbb P^m\) the common vanishing locus of the polynomials in \(W\). Beelen, Datta and Ghorpade consider the quantity
\[
e_r^{\mathbb P}(d,m;q)
:=
\max\bigl\{|V(W)(\mathbb F_q)|:
W\subseteq\mathbb F_q[x_0,\ldots,x_m]_d,\ 
\dim_{\mathbb F_q}W=r\bigr\}.
\]
They conjectured an explicit formula for \(e_r^{\mathbb P}(d,m;q)\) and proved the formula for certain ranges of $r$ in \cite{beelen2018maximum, beelen2022combinatorial}. The problem of determining $e_r^{\mathbb P}(d,m;q)$ is equivalent to computing the $r$-th generalized Hamming weight of the $q$-ary projective Reed–Muller code $\operatorname{PRM}_q(d,m)$.
This problem has been studied extensively, see
\cite{beelen2019vanishing,boguslavsky1997number,
datta2015conjecture,datta2017number,datta2017remarks,
san2024recursive,serre1991lettre,
singhal2025conjecturebeelendattaghorpade,
sorensen1991projective,zanella1998linear}.

When $m=2$, Theorem \ref{thm: cases of urdm} implies the
Beelen--Datta--Ghorpade conjecture for all values of $r$, as shown in
\cite[Section~6]{lin2026largest}. For general $m$, it addresses the
problem of bounding common zero loci that are finite. In coding-theoretic
terms, an upper bound for $u_r(d,m)$ gives a lower bound for the support
size of every $r$-dimensional sub-code of $\operatorname{PRM}_q(d,m)$ whose
defining polynomials have a zero-dimensional common vanishing locus.
Independently of this interpretation, Question~1.1 is a geometric problem
of interest in its own right. It is also the special case $k=0$ in  \cite[Conjecture~1.5]{lin2026k}, which concerns the upper bound for $|V(W)(\Fq)|$ when $\dim(V(W)) \leq k$.

As explained in
\cite[Proposition~3.1]{lin2026k}, Conjecture \ref{Conjecture urdm} follows from the Eisenbud--Green--Harris conjecture introduced in \cite{eisenbud1993authors} through
Artinian reduction. In this paper, we give a direct, elementary proof of Theorem~\ref{Thm: main result}
that does not rely on the E.G.H. conjecture and imposes no restriction
on the characteristic of the base field $\kappa$.
See \cite{caviglia2021eisenbud, gunturkun2021survey} for progress that has been made on the E.G.H. conjecture as well as its applications.

\subsection{Reduction of Theorem~\ref{Thm: main result} to Theorem~\ref{Thm: main result g_Gamma}}
In this subsection we briefly explain how Theorem \ref{Thm: main result} follows from Theorem \ref{Thm: main result g_Gamma}. The argument follows that of \cite[Lemma~5.5, Theorem~5.7]{lin2026largest}.

We first give the expression for $H'_{r-(m-1)}(d,m)$. For $m\leq r\leq\binom{m+d}{d}$ with
$r>\binom{m+d}{d}-\binom{d+2}{2}$, the definition of $H'_{r-(m-1)}(d,m)$ gives
\[
H'_{r-(m-1)}(d,m)=
\begin{cases}
\tbinom{m+d}{d}-r,
& \tbinom{m+d}{d}-3d+1<r\leq\tbinom{m+d}{d},\\[8pt]
\tbinom{m+d}{d}-r+\binom{t-1}{2},
& \tbinom{d-t+1}{2}
  <r-\tbinom{m+d}{d}+\tbinom{d+2}{2}-1
  \leq\tbinom{d-t+2}{2}, \ \ 3\leq t\leq d\\[8pt]
d^2+1,
& r=\tbinom{m+d}{d}-\tbinom{d+2}{2}+1,
\end{cases}
\]

\begin{proof}[Proof of Theorem~\ref{Thm: main result} assuming Theorem~\ref{Thm: main result g_Gamma}]
By Proposition~\ref{Prop: lower bound urdm}, it suffices to prove
$u_r(d,m)\leq H'_{r-(m-1)}(d,m)$.
Suppose that there exists an $r$-dimensional subspace
$W\subseteq S(m)_d$ with $\dim V(W)=0$ and
$|V(W)|>H'_{r-(m-1)}(d,m)$. By \cite[Lemma~5.3]{lin2026largest},
we can choose a sub-scheme $\Gamma\subseteq V(W)$ with
\[
\deg(\Gamma)=H'_{r-(m-1)}(d,m)+1.
\]
By Lemma~\ref{lem: regular sequence Krull dimension}, there exist $G_1,\ldots,G_m\in W$
forming a regular sequence.
Thus $\Gamma\subseteq V(G_1,\ldots,G_m)$.

In the three ranges of $r$ above, respectively, we have
\[
\deg(\Gamma)\leq 3d-1,\qquad
\deg(\Gamma)\leq td+d-t+1,\qquad
\deg(\Gamma)=d^2+2.
\]
The corresponding parts of
Theorem~\ref{Thm: main result g_Gamma} therefore give
\[
g_\Gamma(d)\leq
\begin{cases}
0 & \text{if $\deg(\Gamma) \leq 3d-1$},\\
\binom{t-1}{2} & \text{if $\deg(\Gamma) \leq td+d-t+1$},\\
\binom{d-1}{2} & \text{if $\deg(\Gamma)=d^2+2$}.
\end{cases}
\]
Thus in all 3 ranges of $r$, we get
\[
g_{\Gamma}(d)
\leq H'_{r-(m-1)}(d,m)+r-\tbinom{m+d}{d}
= \deg(\Gamma)-1 +r -\tbinom{m+d}{d}.
\]
Since $W\subseteq I(\Gamma)_d$, we obtain 
\[
r\leq\dim_\kappa I(\Gamma)_d
=\tbinom{m+d}{d}-\deg(\Gamma)+g_\Gamma(d)
\leq r-1,
\] 
which is a contradiction.
\end{proof}
\subsection{Artinian reduction}\label{Subsec: Artinian reduction}
Artinian reduction allows us to express the geometric quantities
$\deg(\Gamma)$ and $g_\Gamma(d)$ in terms of the Hilbert function
of an Artinian $\ka$-algebra $A(\Gamma)$.
We recall the construction described in
\cite{eisenbud1996cayley}, but adapted to our context.

Recall that $S(\Gamma):=S(m)/I(\Gamma)$ is the homogeneous coordinate ring of $\Gamma$ and $S(\Gamma)_k$ is its $k^{th}$ graded piece. 

% Throughout this paper, let $\kappa$ denote an algebraically closed field (possibly of characteristic $0$) and let $S:=\ka[x_0, \dots, x_m]$ be the homogeneous coordinate ring of $\PP^m$. Let $S_d$ be the degree $d$ graded piece of $S$.

% Let $\Gamma \subseteq \PP^m$ be a zero-dimensional scheme.  Let $I(\Gamma) \subseteq S$ be the homogeneous ideal of $\Gamma$, and let $S_{\Gamma}:=S/I(\Gamma)$ be the homogeneous coordinate ring of $\Gamma$. 

Now, let $\ell \in S(m)_1$ be a linear form such that $V(\ell)\cap \Gamma=\emptyset$. Since $\Gamma$ is zero-dimensional, $\sg$ is a one-dimensional Cohen-Macaulay ring. Since $V(\ell)\cap\Gamma=\emptyset$, $\ell$ is not a zero-divisor in $\sg$. The \textbf{Artinian reduction} associated to $\Gamma$ (with respect to $\ell$) is:
\[
A(\Gamma):=S(\Gamma)/(\ell)=S(m)/(I(\Gamma),\ell).
\]
By its construction, $A(\Gamma)$ has Krull dimension $0$, hence it is a local Artinian algebra with the unique maximal ideal $\mathfrak{m}=(x_0, \dots, x_m)$. For a graded algebra $A$, let $A_i$ be its $i$-th graded piece.
We have the following exact sequence:
\[
0 \to \sg_{i-1} \xrightarrow{\cdot \ell}\sg_i \to \ag_{i} \to 0.
\]
Let $h_{\ag}$ (resp. $h_{\sg}$) be the Hilbert function of the corresponding graded algebra. From the above exact sequence, we have
\[
h_{\ag}(i)=h_{\sg}(i)-h_{\sg}(i-1).
\]
Taking the sum over all graded pieces, we get the following properties of $\ag$:
\begin{enumerate}
    \item $\dim_\ka\ag
    =\lim_{N \to \infty}h_{\sg}(N)=\deg(\Gamma)$,
    \item $\dimk {\ag}_{\leq d}
    =h_{\sg}(d)=\binom{m+d}{d}- \dim_{\ka} I(\Gamma)_d$,
    \item $\dimk {\ag}_{\geq d+1}=g_\Gamma(d)$.
\end{enumerate}
For a standard graded Artinian $\kappa$-algebra $A$, we write
$g_A(i):=\dim_\kappa A_{\geq i+1}$ for $i\geq 0$.

Now, suppose $\Gamma\subseteq V(F_1,\ldots,F_m)$, where
$F_1,\ldots,F_m\in S(m)_d$ form a regular sequence.
Choose $\ell\in S(m)_1$ such that
$V(\ell)\cap V(F_1,\ldots,F_m)=\emptyset$ in $\PP^m$.
After a change of coordinates, we may assume that $\ell=x_0$.
Since $S(m)/(x_0,F_1,\ldots,F_m)$ is Artinian,
$x_0,F_1,\ldots,F_m$ form a regular sequence
in the Cohen--Macaulay ring $S(m)$.
Consequently, the images $\Fb_1,\ldots,\Fb_m$
form a regular sequence of degrees $(d,\ldots,d)$ in
$R(m):=S(m)/(x_0)\cong\kappa[x_1,\ldots,x_m]$.

Let $J$ be the image of $I(\Gamma)$ in $R(m)$.
Then $A(\Gamma)=R(m)/J$, and $J$ contains
the regular sequence $\Fb_1,\ldots,\Fb_m$ because
$I(\Gamma)$ contains $F_1,\ldots,F_m$. Hence, to prove Theorem~\ref{Thm: main result g_Gamma}, it suffices to prove the following statements:
\begin{proposition}\label{Prop: algebra reg seq 3d-1}
Let $A=\kappa[x_1,\ldots,x_m]/I$, where $I$ is a homogeneous
ideal containing a regular sequence of degrees $(d,\ldots,d)$.
If $\dim_\kappa A\leq 3d-1$, then $g_A(d)=0$.
\end{proposition}
\begin{theorem}\label{thm: main inequality}
Let $A=\kappa[x_1,\ldots,x_m]/I$, where $I$ is a homogeneous
ideal containing a regular sequence of degrees $(d,\ldots,d)$.
Then, for every $2\leq t\leq d$,
\begin{equation}\label{eq: main inequality}
\dim_\kappa A\leq td+d-t+1
\Longrightarrow
g_A(d)\leq \tbinom{t-1}{2}.
\end{equation}
\end{theorem}
\begin{proposition}\label{prop: end bound}
    Let $A=\kappa[x_1,\ldots,x_m]/I$, where $I$ is a homogeneous
ideal containing a regular sequence of degrees $(d,\ldots,d)$. If $\dimk A \leq d^2+2$, then $g_A(d) \leq \binom{d-1}{2}$.
\end{proposition}
Propositions \ref{Prop: algebra reg seq 3d-1} and \ref{prop: end bound}  will be proven in Section \ref{sec: smallest and largest r}, while Theorem \ref{thm: main inequality} will be proven in Section~\ref{Section: Artinian tail}.

\subsection{Lemmas on Artinian algebra}\label{Section: hilbert function}
In this sub-section, we recall some lemmas on Artinian algebras. Most of these results have more general formulations. Here, we state the corollaries that we will use in this paper and provide self-contained proofs.

We first recall the following criteria for extracting a regular sequence of polynomials from a $\kappa$ vector subspace $W$ contained in some graded $\kappa$-algebra. We will use this lemma extensively in this paper.
\begin{lemma}\cite[Lemma 2.1]{lin2026k}\label{lem: regular sequence Krull dimension}
Let $A$ be a graded Cohen-Macaulay $\ka$-algebra of Krull dimension $m$. 
Let $W\subseteq A_d$ be a
finite-dimensional vector subspace for some $d\geq 1$, and let $I(W)$ be the ideal generated by $W$. Let $B=A/I(W)$. If the Krull dimension of $B$ is at most $m-k$, then we can find a regular sequence $(F_1,F_2, \dots, F_{k})$ in~$W$.
\end{lemma}
\subsubsection{Gorenstein quotients}
We introduce the concepts of socle degree and Gorenstein quotient, following \cite[Section 1.3]{eisenbud1996cayley}.

Recall that the Artinian reduction $A=A(\Gamma)$ is a local
Artinian $\kappa$-algebra with maximal ideal
$\mathfrak m=(x_0,\ldots,x_m)$ and residue field
$A/\mathfrak m\cong\kappa$. For such local Artinian algebras $A$, 
we call the largest integer $k_0$ such that $A_{k_0}\neq 0$
the \textbf{socle degree} of $A$. We say that
$A$ is \textbf{Gorenstein} if
$\operatorname{Ann}_A(\mathfrak m)$ is one-dimensional over
$\kappa$.

% Following \cite[Section~1.3]{eisenbud1996cayley}, we say that
% $A$ is \textbf{Gorenstein} if
% $\operatorname{Ann}_A(\mathfrak m)$ is one-dimensional over
% $\kappa$; in this case, the socle is precisely $A_{k_0}$.

If $\Gamma$ is a zero-dimensional complete intersection, then
its Artinian reduction $A(\Gamma)$ is Gorenstein, although
the converse does not hold in general.
Among the equivalent characterizations given in
\cite[Section~1.3]{eisenbud1996cayley}, we will use the
following property.

\begin{lemma}\cite[Proposition 7]{eisenbud1996cayley}\label{lem:gorenstein-symmetry}
Let $G$ be a standard graded Artinian Gorenstein algebra of
socle degree $k_0$. Then $G_{k_0}$ is one-dimensional, and
for every $0\leq i\leq k_0$, multiplication induces a perfect
pairing
\[
G_i\times G_{k_0-i}\longrightarrow G_{k_0}\cong\kappa.
\]
Consequently, $h_G(i)=h_G(k_0-i)$ for every $0\leq i\leq k_0$,
and $h_G(0)=h_G(k_0)=1$.
\end{lemma}

We will also use the fact that every standard graded Artinian algebra has a Gorenstein quotient
with the same socle degree.

\begin{lemma}\cite[Proof of Proposition 12]{eisenbud1996cayley}\label{lem: gorenstein-quotient}
Let $A$ be a standard graded Artinian algebra of socle degree
$k_0$. Then $A$ admits a graded Gorenstein quotient $G$ of the same
socle degree $k_0$.
\end{lemma}

% If $A$ is the Artinian reduction of a zero-dimensional scheme $\Gamma$, one can think of $G$ arised from the smallest complete intersection in $\PP^m$ that contains $\Gamma$.

\subsubsection{Macaulay's growth condition}

We follow the exposition in
\cite[Section~2]{gunturkun2021survey}.

We first define
the Macaulay transform.
For positive integers $n$ and $i$, write the unique
$i$-binomial expansion
\[
n=\tbinom{n_i}{i}+\tbinom{n_{i-1}}{i-1}
+\cdots+\tbinom{n_j}{j},
\]
where $n_i>n_{i-1}>\cdots>n_j\geq j\geq 1$.
Define the degree-$i$ Macaulay transform of $n$ by
\[
n^{\langle i\rangle}
=
\tbinom{n_i+1}{i+1}
+\tbinom{n_{i-1}+1}{i}
+\cdots+\tbinom{n_j+1}{j+1}.
\]
We also set $0^{\langle i\rangle}=0$.

From the definition of $n^\itr$ it is not hard to derive the following lemma.
\begin{lemma}\label{lem: macualy transform}
Let $s \in \N$. We have the following properties of the $i^{th}$ Macaulay transform of $s$:
\[
s^{\langle i \rangle} = \begin{cases}
s & \text{if } s\leq i,\\
s+1 & \text{if } i+1 \leq s\leq 2i.
\end{cases}
\]
\end{lemma}

Macaulay's growth condition provides an upper bound for the Hilbert function of a standard graded algebra at each degree, depending on the value at the preceding degree.
\begin{theorem}[Macaulay's growth condition]\cite[Theorem 2.6]{gunturkun2021survey}
\label{thm:macaulay}

\noindent Let $A$ be a standard graded algebra. Then, for every
$i\geq 1$, we have $h_A(i+1)\leq h_A(i)^{\langle i\rangle}$.
In particular, if $h_A(i)\leq i$, then $h_A(i+1)\leq h_A(i)$.
\end{theorem}

We will use the following corollaries of Macaulay's growth condition:
\begin{corollary}\label{cor: growth condition cor}
Let $A$ be a standard graded Artinian algebra. Then, for $i\leq j$, we have
\begin{enumerate}
    \item $h_A(i) \geq \min(i+1, h_A(j))$;
    \item $h_A(i) \geq \min(2i+1, h_A(j) -(j-i))$.
\end{enumerate}
\end{corollary}
\begin{proof}
\noindent \begin{enumerate}
\item By Lemma \ref{lem: macualy transform}, if $h_A(i) \le i$, then $h_A(i+1) \le h_A(i)^{\langle i \rangle} = h_A(i)$, and so $h_A(j) \le h_A(i)$ for all $j \ge i$. This shows $(1)$.
\item If $h_A(i) \leq i$ then the statement follows from $(1)$. If $i+1 \leq h_A(i) \leq 2i$, by Lemma \ref{lem: macualy transform},  $h_A(i+1) \leq h_A(i)^{\itr}=h_A(i)+1$. It follows that $h_A(j) \leq h_A(i)+j-i$. This shows $(2)$.
\end{enumerate}
\end{proof}

\begin{corollary}\label{cor: macaulay lower bound}
Let $B$ be a standard graded Artinian algebra over $\ka$. Let $h_B(i)$ be the Hilbert function of $B$ at degree $i$. 
Suppose $h_B(d)=s$. Then,
\[
\dim B_{\le d-1} \geq
\begin{cases}
sd-\binom{s}{2} 
& \text{if } 1 \leq s\leq 2d;\\
\binom{d+1}{2}+(s-d-1)(d-1)-\binom{s-d-1}{2} 
&\text{if } d+1 \leq s\leq 2d;\\
d^2
&\text{if } s\geq 2d+1.
\end{cases}
\]
\end{corollary}

\begin{proof}
For $1\le s\le d+1$, Corollary~\ref{cor: growth condition cor}(1) gives
\[
\dim_\kappa B_{\le d-1}
\ge \tsum_{i=0}^{d-1}\min\{i+1,s\}
=sd-\tbinom{s}{2}.
\]
For $s\ge d+1$, Corollary~2.9(2) gives, for
$0\le i\le d-1$,
\[
h_B(i)\ge \min\{2i+1,s-d+i\}
=i+1+\min\{i,s-d-1\}.
\]
Summing these inequalities yields
\[
\begin{aligned}
\dim_\kappa B_{\le d-1}
&\ge \tsum_{i=0}^{d-1}
   \bigl(i+1+\min\{i,s-d-1\}\bigr)\\
&=
\begin{cases}
\displaystyle
\tbinom{d+1}{2}+(s-d-1)(d-1)-\tbinom{s-d-1}{2},
& d+1\le s\le 2d,\\[4pt]
d^2, & s\ge 2d+1.
\end{cases}
\end{aligned}
\]
Finally, for $d+1\le s\le 2d$, the middle bound also
implies the first bound, since
\[
\tbinom{d+1}{2}+(s-d-1)(d-1)-\tbinom{s-d-1}{2}
-\left(sd-\tbinom{s}{2}\right)
=d(s-d-1)\ge0.
\]
\end{proof}

We next recall Gotzmann's persistence theorem. It basically says that if the equality in Macaulay's growth condition holds at a certain degree $i$, then the equality holds for all the higher degrees. 
\begin{theorem}[Gotzmann persistence]
\label{thm:gotzmann}\cite[Theorem 2.8]{gunturkun2021survey}

\noindent Let $A=\kappa[x_1,\ldots,x_m]/I$, where $I$ is a homogeneous
ideal generated by elements of degree at most $r+1$, for some integer
$r\geq 1$. If $h_A(r+1)=h_A(r)^{\langle r\rangle}$, then $h_A(i+1)=h_A(i)^{\langle i\rangle}$ for every $i\geq r$.

In particular, if $A$ is Artinian and $h_A(i)>0$, then we must have $h_A(i+1)<h_A(i)^\itr$ for $i \geq r$.
\end{theorem}

We will use the following simple corollary throughout the paper.
\begin{corollary}\label{cor: Gotzman persistence}
Let $A=\kappa[x_1,\dots,x_m]/I$ be a standard graded Artinian algebra, where $I$ contains a regular sequence of degree $(d,\dots,d)$ and $d\ge2$. If $i\ge d-1$ and $h_A(i)>0$, then $h_A(i+1)<h_A(i)^{\langle i\rangle}$. In particular, if $0<h_A(i)\le i$, then $h_A(i+1)<h_A(i)$.
\end{corollary}

\begin{proof}
By Macaulay's growth condition, it suffices to rule out equality. Suppose that $h_A(i+1)=h_A(i)^{\langle i\rangle}$. Let $J$ be the ideal generated by the homogeneous components of $I$ of degrees at most $i$, and set $A'=\kappa[x_1,\dots,x_m]/J$. Since $J\subseteq I$ and $J_i=I_i$, we have $h_{A'}(i)=h_A(i)$. Applying Macaulay's growth condition to $A'$ gives
\[
h_A(i+1)\le h_{A'}(i+1)
\le h_{A'}(i)^{\langle i\rangle}
=h_A(i)^{\langle i\rangle}
=h_A(i+1).
\]
Thus equality holds throughout, so $J_{i+1}=I_{i+1}$. Since $J$ and $I$ already agree in degrees at most $i$ and $d\le i+1$, we obtain $J_d=I_d$. Hence $J$ contains the given regular sequence, and $A'$ is Artinian.

On the other hand, $J$ is generated in degrees at most $i$ and $h_{A'}(i+1)=h_{A'}(i)^{\langle i\rangle}$. Gotzmann persistence therefore gives $h_{A'}(j+1)=h_{A'}(j)^{\langle j\rangle}$ for every $j\ge i$. Since $h_{A'}(i)>0$, these values remain positive, contradicting the fact that $A'$ is Artinian.

Finally, if $0<h_A(i)\le i$, then $h_A(i)^{\langle i\rangle}=h_A(i)$, giving the last assertion.
\end{proof}

% \section{Artinian algebra result implies geometric results}\label{Section: algebra implies geometry}

% \begin{corollary}\label{cor: no d+2 in a line and no 2d+2 on a curve}
% Let $\Gamma \subseteq V(F_1, \dots, F_m) \subseteq \PP^m$ be a subscheme of complete intersection of degree $d$ hypersurfaces, such that $\deg(\Gamma) \leq 2d+2$. Then $g_{\Gamma}(d)=0$.
% \end{corollary}
% \begin{proof}
% We just need to verify that $\Gamma$ does not have $d+2$ points in a line and $2d+2$ points on a conic.

% If $\Gamma$ has $d+1$ points on a line $V(L)$, then the $d+1$ points lie in $V(L) \cap V(F_i)$ for each $F_i$. 
% If $V(L)$ is not contained in $V(F_i)$, then $V(L) \cap V(F_i)$ has at most $d$ points. 
% Therefore, $V(L) \subseteq V(F_i)$ for all $1 \leq i \leq m$. 
% Hence $L \subseteq  V(F_1, \dots, F_m)$, this contradicts the fact that $(F_1, \dots, F_m)$ is a regular sequence.

% Similarly, if $\Gamma$ has $2d+1$ points on a conic $V(Q)$, then by the same reasoning $V(Q) \subseteq V(F_1, \dots, F_m)$, which contradicts that $(F_1, \dots, F_m)$ is a complete intersection.
% \end{proof}

% \begin{theorem}\label{thm: main inequality}
% Let $(F_1, \dots, F_m) \subseteq S(m)_d$ be a regular sequence of degree $(d, \dots, d)$ and let $\Gamma \subseteq V(F_1, \dots, F_m)$ be a subscheme. Then, for every $3 \leq t \leq d$,
% \begin{equation}\label{eq: main inequality}
% \deg(\Gamma)\le td+d-t+1 \Longrightarrow g_\Gamma(d)\le \binom{t-1}{2}.
% \end{equation}
% \end{theorem}

\section{Reduction from Theorem \ref{Thm: 3d-1 geometric} to Theorem \ref{thm: 3d-1 artinian}}\label{sec: reduction from 1.8 to 3.1}

To prove Theorem \ref{Thm: 3d-1 geometric}, we reduce it to the following statement about graded quotients of Artinian algebras.
\begin{theorem}\label{thm: 3d-1 artinian}
Let $A$ be a standard graded Artinian algebra such that $\dimk A\leq 3d-1$. 
Then $g_A(d)=0$ if and only if
every graded quotient $B$ of $A$ satisfies both of the following conditions:
\begin{itemize}
    \item if $h_B(1)\leq 1$, then $\dimk B \leq d+1$;
    \item if $h_B(1)=2$ and $h_B(2) \leq 2$, then $\dimk B \leq 2d+1$.
\end{itemize}
\end{theorem}

Theorem \ref{thm: 3d-1 artinian} will be proven in Section~\ref{sec: smallest and largest r}. Here, we show that it implies Theorem \ref{Thm: 3d-1 geometric}.

\begin{proof}[Proof of Theorem \ref{Thm: 3d-1 geometric} assuming Theorem \ref{thm: 3d-1 artinian}]

First, consider a zero-dimensional scheme $\Gamma \subseteq \PP^m$ for which $g_\Gamma(d)=0$ and $\deg(\Gamma) \leq 3d-1$. Choose $\ell \in \sg_1$ such that $V(\ell) \cap \Supp(\Gamma)=\emptyset$ and let $\ag=\sg/\ell$ be the Artinian reduction of $\Gamma$. Then $g_{\ag}(d)=g_\Gamma(d)=0$. 

\begin{itemize}
    \item Given a line $L\subseteq \Pm$, let $I(L)\subseteq S(m)$ be its homogeneous ideal. Then, $A(\Gamma \cap L)$ is a graded quotient of $\ag$ that factors through $\ag/\overline{I(L)}$. Hence, $h_{A(\Gamma \cap L)}(1) \leq 1$. 
    By Theorem~\ref{thm: 3d-1 artinian}, we have $\dimk A(\Gamma \cap L)\leq d+1$, so $\deg(\Gamma \cap L) \leq d+1$.
    \item Suppose $C$ is a planar conic with homogeneous ideal $I(C)$ in $S(m)$.
    For this case, choose $\ell\in S(m)_1$ such that $V(\ell)\cap \Supp(\Gamma)=\emptyset$ and $V(\ell)$ contains no irreducible component of $C$.
    Form both $A(\Gamma)$ and $A(\Gamma\cap C)$ using this same $\ell$.
    We know that $\dimk(I(C)_1)=m-2$ and there exists $Q \in I(C)_2$ such that $Q\notin \langle I(C)_1, \ell\rangle$. Then, $A(\Gamma \cap C)$ is a graded quotient of $\ag$ that factors through $\ag/(\overline{I(C)})$. Hence, $h_{A(\Gamma \cap C)}(1) \leq 2$ and $h_{A(\Gamma \cap C)}(2) \leq 2$. 
    If $h_{A(\Gamma \cap C)}(1) \leq 1$ then 
    \[
    \deg(\Gamma \cap C) =\dimk A(\Gamma\cap C) \leq d+1 \leq 2d+1.
    \] 
    If $h_{A(\Gamma \cap C)}(1) =2$, then $\deg(\Gamma \cap C) = \dimk A(\Gamma\cap C) \leq 2d+1$.
\end{itemize}
This completes the proof of one direction.

Now, suppose $\Gamma \subseteq \PP^m$ is a zero-dimensional scheme such that $\deg(\Gamma)\leq 3d-1$ and $g_\Gamma(d)>0$. 
We want to find a line $L$ or a conic $C$ whose intersection with $\Gamma$ exceeds the corresponding bound.

Without loss of generality, we replace $\Gamma$ by a subscheme of $\Gamma$ of minimal degree that fails to impose independent conditions on polynomials of degree $d$.
By abuse of notation, we continue to denote this sub-scheme by $\Gamma$. By minimality, for any proper sub-scheme $\Gamma' \subsetneq \Gamma$ such that $\deg(\Gamma')=\deg(\Gamma)-1$, $g_{\Gamma'}(d)=0$. Hence, we have $g_\Gamma (d)=1$. We will show that this minimal length sub-scheme $\Gamma$ is contained either in a line $L$ or a conic $C$.

Let $\sg$ be the coordinate ring of $\Gamma$ and choose $\ell \in S(m)_1$ such that $\Gamma \cap V(\ell) =\emptyset$. Let $\Ag:=\sg/\ell \sg$ be the Artinian reduction. Since $\ga(d)=1$, we must have $h_{\ag}(d+1)=1$ and $h_{\ag}(i)=0$ for $i \geq d+2$. Consequently, for any $r \geq 1$, the map given by multiplication by $\ell^r$, $\ell^r: \sg_{i} \to \sg_{i+r}$
is an isomorphism for $i \geq d+1$. 

We first establish the following two claims on $\sg$ and $\ag$.

\noindent\textbf{Claim 1}:
Let $e \geq 1$, and $f \in \sg_e$. If $f\cdot\sg_{d+1} \subseteq \ell^{e+1}\sg_{d}$, then $f=0$.\\
\noindent\textbf{Proof of Claim 1}:
Assume for the sake of contradiction that $f \neq 0$. Set $\Gamma_1=V(f) \cap \Gamma$. Since $f \neq 0$, we have $\Gamma_1 \subsetneq \Gamma$. We will show that $\ggo(d)=1$, contradicting the minimality of $\Gamma$.

Consider the maps on the coordinate rings induced by restriction to $\go$:
\begin{align*}
&\Res_{d+1+e}: \sg_{d+1+e} \twoheadrightarrow \sgo_{d+1+e},
&
&\Res_{d}: \sg_{d} \twoheadrightarrow \sgo_{d}.
\end{align*}
Note that $\dimk(\sg_{d+1+e}) - \dimk(\sg_d) = g_{\Gamma}(d)=1$ and 
\begin{align*}
g_{\Gamma_1}(d)
&\geq \dimk(\sgo_{d+1+e}) - \dimk(\sgo_d)\\
&= \dimk(\sg_{d+1+e}) -\dimk(\ker(\Res_{d+1+e}))
-\dimk(\sg_{d}) +\dimk(\ker(\Res_{d}))\\
&=1 -\dimk(\ker(\Res_{d+1+e}))+\dimk(\ker(\Res_{d})).
\end{align*}
We will show that $\ker(\Res_{d+1+e})=f\cdot\sg_{d+1}$. Indeed, $\sgo=S(m)/(I(\Gamma),f)^{\text{sat}}$, so $\ker(\Res_{d+1+e}) \supseteq f\cdot\sg_{d+1}$. 
Conversely, let $u \in \sg_{d+1+e}$ such that $u \in \overline{(I(\Gamma),f)^{\text{sat}}}$. 
Then there exists $r \geq 0$ such that $u\ell^r=fv$ for some $v \in \sg_{d+r+1}$. Since $\ell^r: \sg_{d+1} \to \sg_{d+1+r}$ is surjective, $v= \ell^r v_1$ for some $v_1\in \sg_{d+1}$. Moreover $\ell$ is a non-zero divisor in $\sg$, so $u=f v_1 \in f\cdot\sg_{d+1}$.

This implies that $\ker(\Res_{d+1+e}) \subseteq \ell^{e+1} \sg_d$. If we restrict the map $\Res_{d+1+e}$ to $\ell^{e+1} \sg_d$, then the kernel does not change, and the image becomes $\ell^{e+1} \sgo_{d}$.
We get the following commutative diagram.
\[
\begin{tikzcd}[row sep=2em, column sep=8em]
\sg_d
    \arrow[r, two heads, "\Res_d"]
    \arrow[d, "\cdot\,\ell^{e+1}"', "\cong"]
&
\sgo_d
    \arrow[d, "\cdot\,\ell^{e+1}", "\cong"']
\\
\ell^{e+1}\sg_d
    \arrow[r, two heads,
        "\left.\Res_{d+1+e}\right|_{\ell^{e+1}\sg_d}"']
&
\ell^{e+1}\sgo_d
\end{tikzcd}
\]
Thus $\dimk(\ker(\Res_d)) = \dimk(\ker(\restr{\Res_{d+1+e}}{\ell^{e+1}\sg_d})) = \dimk(\ker(\Res_{d+1+e}))$. Hence, $g_{\Gamma_1}(d)\geq 1$, which is a contradiction. We conclude that $f=0$, and the proof of claim 1 is complete.

% Now, restrict the map $\Res_{\go}$ to the subspace $\ell^{e+1}\sg_{d}$. Then, $\ker(\Res_{\go}) \cap \ell^{e+1}\sg_{d}=f(\sg)_{d+1}$, and $\Res_{\go}(\ell^{e+1}\sg_d)=\ell^{e+1}(\sgo)_d$. Since multiplication by $\ell$ is injective on each graded piece of $\sg$ and $\sgo$, we have
% \begin{align*}
%     h_{\sgo}(d)&=\dim((\sgo)_d)=\dim(\ell^{e+1}(\sgo)_d)=\dim(\ell^{e+1}\sg_d)-\dim(f(\sg)_{d+1})\\
%     &=h_{\sg}(d)-\dim(f(\sg)_{d+1})\\
%     h_{\sgo}(d+e+1)&=h_{\sg}(d+e+1)-\dim(f(\sg)_{d+1})
% \end{align*}
% Subtracting the two equations, we have $h_{\sgo}(d+e+1)-h_{\sgo}(d)=h_{\sg}(d+e+1)-h_{\sg}(d)$. Since $\Ag_{\geq d+2}=0$, $(A_{\Gamma_1})_{\geq d+2}=0$, so for any $i \geq d+2$, $h_{\sgo}(i)=\deg(\Gamma_1)$. In particular, $h_{\sg}(d+e+1)=\deg(\Gamma)$. Therefore, \[\ggo(d)=h_{\sgo}(d+e+1)-h_{\sgo}(d)=h_{\sg}(d+e+1)-h_{\sg}(d)=g_{\Gamma}(d)=1\] This contradicts the minimality of $\Gamma$, so $f=0$. Hence, the proof of claim 1 is complete.

\noindent\textbf{Claim 2}:
Fix $0\leq k\leq d+1$. If for all $0 \leq i <k$, the homomorphism
\begin{align*}
    \Ag_i &\to \text{Hom}_\ka(\Ag_{d+1-i}, \Ag_{d+1})\\
    a_i & \mapsto (a \mapsto aa_i)
\end{align*}
is surjective, then for any standard graded quotient $q: \ag \to B$ such that $h_B(d+1) \neq 0$, $q:\Ag_k \to B_k$ is an isomorphism.

\noindent\textbf{Proof of Claim 2}:
The case $k=0$ is immediate, so assume $k\geq1$.
It suffices to show that $q:\Ag_k \to B_k$ is injective. Suppose $a \in \Ag_k$ is such that $q(a)=0$. Then, for every $a' \in \Ag_{d+1-k}$, we have $q(a a')=0$ in $B_{d+1}$. Since $h_B(d+1) \neq 0$, $q: A(\Gamma)_{d+1} \to B_{d+1}$ is an isomorphism. Hence, $a\Ag_{d+1-k}=0$. 

We will construct a lift of $a$ to $S(\Gamma)_k$ and then use Claim 1 to show that this lift vanishes and hence $a=0$.
For $1 \leq i \leq k+1$, we inductively construct lifts $f_i \in \sg_k$ such that $f_i \cdot\sg_{d+i-k} \subseteq \ell^i \sg_d$.

Start with an arbitrary lift $f \in \sg_k$ of $a$. Since $a\Ag_{d+1-k}=0$, we see that $f \cdot \sg_{d+1-k} \subseteq \ell\sg_d$. Hence, we can take $f_1=f$ and this finishes the base case.

Now suppose that for some $1 \leq i \leq k$, we have the lift $f_i$ such that $f_i\cdot \sg_{d+i-k} \subseteq \ell^{i} \sg_{d}$. Since multiplication by $\ell^{i}$ is an isomorphism $\sg_{d+1} \to \sg_{d+1+i}$, consider $\ell^{-i}f_i: \sg_{d+i+1-k} \to \sg_{d+1}$. On the subspace $\ell \sg_{d+i-k}$, the image of $\ell^{-i}f_i$ is in $\ell \sg_d$. Hence, this map descends to $[\ell^{-i}f_i]:\Ag_{d+i+1-k} \to \Ag_{d+1}$. By our assumption, there exists $c_{k-i} \in \Ag_{k-i}$ such that $[\ell^{-i}f_i]=c_{k-i}$. That is, $[\ell^{-i}f_i-c_{k-i}]$ is zero on $\Ag_{d+i+1-k}$. Let $\Tilde{c}_{k-i}$ be a lift of $c_{k-i}$ in $\sg_{k-i}$. Then $\ell^{-i}f_i-\Tilde{c}_{k-i}$ maps $\sg_{d+i+1-k}$ to $\ell\sg_d$. Set $f_{i+1}=f_i-\tilde{c}_{k-i}\ell^{i}$, then $f_{i+1} \cdot \sg_{d+i+1-k}\subseteq \ell^{i+1}\sg_d$. This establishes the inductive step.

Now let $i=k+1$. We have that $f_{k+1} \cdot \sg_{d+1} \subseteq \ell^{k+1}\sg_d$. But then by Claim 1, $f_{k+1}=0$. Since $f_{k+1}$ is also a lift of $a$, we have $a=0$. Hence, $q:A(\Gamma)_k \to B_k$ is injective and this completes the proof of claim 2. 

We now apply these claims to the quotient $B$ given by Theorem \ref{thm: 3d-1 artinian}. Since $g_{\ag}(d)=1>0$, by Theorem \ref{thm: 3d-1 artinian}, there exists a graded quotient $q: \Ag \to B$ such that one of the following is true:
\begin{enumerate}
    \item $h_B(1)\leq 1$, $\dimk B \geq d+2$;
    \item $h_B(1)=2, h_B(2) \leq 2$ and $\dimk B \geq 2d+2$.
\end{enumerate}
\noindent\textbf{Claim 3}: In both cases $(1)$ and $(2)$, $h_B(d+1)=1$.

\noindent\textbf{Proof of Claim 3}:
    Since $B$ is a quotient of $\ag$, it suffices to show that $h_B(d+1) \neq 0$.
    
    If $h_B(1) \leq 1$, by Corollary \ref{cor: growth condition cor}, we have $h_B(j) \leq h_B(1) \leq 1$ for $1 \leq j$. If $h_B(d+1)=0$, then $\dimk B \leq d+1<d+2$.

    If $h_B(1)=2$ and $h_B(2) \leq 2$, by Corollary \ref{cor: growth condition cor}, we have $h_B(j) \leq h_B(2) \leq 2$ for $2 \leq j$. If $h_B(d+1)=0$, then $\dimk B \leq 2d+1<2d+2$. This completes the proof of Claim 3.

We now apply Claim 1-3 in cases $(1)$ and $(2)$.

\noindent\textbf{Case $(1)$}: Let $q: \ag \to B$ be such a quotient. Then $h_B(1)=1$. Let $k=1$ and apply Claim 2. The assumption is automatic, so we get that $q: \Ag_1 \to B_1$ is an isomorphism. Consequently, $h_{\ag}(1)=1$, so $h_{\sg}(1)=2$. Therefore, $\dimk I(\Gamma)_1=m-1$, which means $\Gamma$ is contained in a line~$L$ and 
\[
\deg(\Gamma\cap L)=\deg(\Gamma)=\dimk A(\Gamma)\geq \dimk B\geq d+2.
\]

\noindent\textbf{Case $(2)$}: 
Let $q: \ag \to B$ be such a quotient. Then $h_B(1)=2$ and $h_B(2) \leq 2$. Let $k=1$ and apply Claim 2. We get that $h_{\ag}(1)=2$. 

Let $k=2$. To apply Claim 2, we need to show that $A(\Gamma)_1 \to \text{Hom}_\ka(\Ag_{d}, \ag_{d+1})$ is surjective. Suppose $a \in A(\Gamma)_1$ is such that $aA(\Gamma)_d=0$, then by the argument of Claim 2, we have $a=0$. Hence, the map $\ag_1 \to \text{Hom}_\ka(\Ag_{d}, \Ag_{d+1})$ is injective. To show that it is surjective, it suffices to show that $h_{\ag}(d)  \leq h_{\ag}(1)$. Indeed, if $h_{\ag}(d) \geq 3$, then by Corollary \ref{cor: growth condition cor}, for $2 \leq i \leq d$, we have $h_{\ag}(i) \geq \min(i+1, h_{\ag}(d)) \geq 3$. Then, $\dimk \ag \ge 1+2+3(d-1) =3d >3d-1$, contradicting our assumption that $\deg(\Gamma) \leq 3d-1$. Hence, $h_{\ag}(d)  \leq 2=h_{\ag}(1)$ and the map is surjective. Applying Claim 2 with $k=2$, we see that $q: \Ag_2 \to B_2$ is an isomorphism, so in particular $h_{\ag}(2)=h_B(2) \leq 2$. 

Therefore, $h_{\sg}(1)=3$ and $h_{\sg}(2) \leq 5$. That means $\dimk I(\Gamma)_1 =m-2$, so $\Gamma$ is contained in a plane $P$. In addition, there is a degree $2$ polynomial $Q \in I(\Gamma)_2$ that is not zero when restricted to the plane $P$, so $\Gamma$ is contained in the plane conic $C=P \cap V(Q)$.

We have therefore shown that the chosen minimal subscheme
$\Gamma$ is either contained in a line $L$ with
$\deg(\Gamma)\geq d+2$, or contained in a planar conic $C$
with $\deg(\Gamma)\geq 2d+2$.
Consequently, the original scheme also has an intersection with a line or a planar conic exceeding the corresponding bound. This proves the converse and hence finishes the proof of Theorem \ref{Thm: 3d-1 geometric}.
\end{proof}

\section{Criterion for $g_A(d)=0$}\label{sec: smallest and largest r}

In this section, we prove Theorem \ref{thm: 3d-1 artinian}, which by the previous section implies Theorem \ref{Thm: 3d-1 geometric}. As corollaries, we deduce Propositions \ref{Prop: algebra reg seq 3d-1} and \ref{prop: end bound}, which yields Part $(1), (3)$ of Theorem \ref{Thm: main result g_Gamma}.
% \begin{customthm}{\ref{thm: 3d-1 artinian}}
%     Let $A$ be a standard graded Artinian algebra over a field $\kappa$ such that $\dim_\kappa A\leq 3d-1$. The following are equivalent:
% \begin{enumerate}
%     \item $g_A(d)=0$.
%     \item Every graded quotient $B$ of $A$ satisfies both
%     of the following conditions:
%     \begin{itemize}
%         \item If $h_B(1)\leq 1$, then $\dim_{\kappa}B \leq d+1$.
%         \item If $h_B(1)=2$ and $h_B(2) \leq 2$, then $\dim_{\kappa}B \leq 2d+1$
%     \end{itemize}
% \end{enumerate}
% \end{customthm}

\begin{proof}[Proof of Theorem \ref{thm: 3d-1 artinian}]

\noindent We first prove that the two quotient conditions imply $g_A(d)=0$. Let $A$ be a standard graded Artinian algebra such that $\dim_{\kappa}A \leq 3d-1$ and every graded quotient $B$ satisfies the conditions. Let $k_0$ be the socle degree of $A$. We want to show that $k_0 \leq d$. 

By Lemma \ref{lem: gorenstein-quotient}, we can form a Gorenstein quotient $G$ of $A$ with the same socle degree $k_0$. Then $\dim_{\kappa}G \leq 3d-1$ and $G$ satisfies the two conditions in Theorem \ref{thm: 3d-1 artinian}. We separate cases based on $h_G(1)$.
\begin{enumerate}
    \item If $h_G(1)=0$, then $G=\kappa$ and $k_0=0\leq d$.
    \item Suppose $h_G(1) \geq 3$. Then, since $G$ is Gorenstein, $h_G(k_0-1) \geq 3$. By Corollary \ref{cor: growth condition cor}, we have that $h_G(i) \geq 3$ for every $2 \leq i \leq k_0-1$. Hence, $\dimk G \geq 2+3(k_0-1)=3k_0-1$. On the other hand, we know that $\dimk G \leq 3d-1$. Hence, $3k_0-1 \leq 3d-1$, so $k_0 \leq d$.
    \item Suppose $h_G(1)=1$. Then, $G \cong \kappa[x]/x^t$ for some $t$, so $k_0=t-1$ and $\dimk G=t$. By assumption on $G$, $\dimk G \leq d+1$, so $k_0 \leq d$. 
    \item Suppose $h_G(1)=2$. Then $h_G(k_0-1)=2$. We separate two cases based on $h_G(2)$.
    \begin{enumerate}
        \item Suppose $h_G(2) \leq 2$. Then, by assumption on $G$, we see that $\dimk G \leq 2d+1$. On the other hand, since  $h_G(k_0-1)=2$, by Corollary \ref{cor: growth condition cor}, $h_G(i) \geq 2$ for $1 \leq i \leq k_0-1$. Therefore, $\dimk G \geq 2(k_0-1)+1+1=2k_0$. Therefore, $2k_0 \leq 2d+1$, so $k_0 \leq d$.
        \item Suppose $h_G(2) \geq 3$. Then, $h_G(k_0-2) \geq 3$. By Corollary \ref{cor: growth condition cor}, $h_G(i) \geq 3$ for $2 \leq i \leq k_0-2$. So, $\dimk G \geq 3(k_0-3)+4+2 =3k_0-3$. Therefore, $3k_0 -3 \leq 3d-1$, so $k_0 \leq d$.
    \end{enumerate}
\end{enumerate} This finishes the proof of $k_0 \leq d$.

Conversely, suppose that $g_A(d)=0$, and let $B$ be a graded quotient of $A$. Let $k_B$ be the socle degree of $B$. Then $k_B \leq d$.
\begin{enumerate}
    \item We show that if $h_B(1) \leq 1$ then $\dimk B \leq d+1$. Since $h_B(1) \leq 1$, we have $B \cong \kappa[x]/x^t$ for some $t$ and $k_B=t-1$. Since $k_B \leq d$, we have $\dimk B=t=k_B+1 \leq d+1$.
    \item We show that if $h_B(1)=2$ and $h_B(2) \leq 2$ then $\dimk B \leq 2d+1$.
    By Corollary~\ref{cor: growth condition cor} with $i=2$, we have $2\geq h_B(2)\geq\min\{3,h_B(j)\}$ for every $j\geq2$.
    Thus $h_B(j)\leq2$ for all $j\geq1$, and $\dim_\kappa B\leq1+2k_B\leq1+2d$.
\end{enumerate}
This finishes the proof.
\end{proof}

We now impose the condition that the defining ideal of the Artinian algebra contains a regular sequence of degree $(\underbrace{d, \dots, d}_m)$. The following lemma will be used in the proofs of both propositions.

\begin{lemma}\label{lem: artinian ring intersection bound}
Let $\kappa$ be an infinite field, and let
$A=\kappa[x_1,\ldots,x_m]/I$, where $I$ is a homogeneous ideal containing a regular sequence
$F_1,\ldots,F_m$ of degrees $(d,\ldots,d)$, with $d\geq 1$.
The following hold for every graded quotient $B$ of $A$:
\begin{enumerate}
    \item If $h_B(1)\leq 1$, then $\dim_\kappa B\leq d$.
    \item If $h_B(1)=2$ and $h_B(j)\leq j$ for some integer
    $j\geq 1$, then $\dim_\kappa B\leq jd$.
\end{enumerate}
In particular, if $h_B(1)=2$ and $h_B(2)\leq 2$, then
$\dimk B\leq 2d$.
\end{lemma}

\begin{proof}
Set $R=\kappa[x_1,\ldots,x_m]$ and let $B$ be a graded quotient of $A$. Then $B=R/I_B$ where $I_B \supseteq  I \supseteq (F_1, \dots, F_m)$. Let $(I_B)_i$ denote the vector space of degree $i$ elements in $I_B$.
\begin{enumerate}
    \item Suppose $h_B(1) \leq 1$. If $h_B(1)=0$, then $\dimk B=1$ and the result follows. Hence, we can assume that $h_B(1)=1$. Choose a basis $\langle\ell_2, \ell_3, \dots, \ell_m\rangle$ of $(I_B)_1$ and extend it to a basis $\langle \ell_1, \ell_2, \dots, \ell_m\rangle$ of $R_1$. Consider the quotient $R/(\ell_2, \dots, \ell_m) \cong \kappa[\ell_1]$. It is Cohen-Macaulay of dimension $1$. Let $\Wb$ be the image of the vector space $\langle F_1, \dots, F_m\rangle$ in $\kappa[\ell_1]$. Then the quotient $\kappa[\ell_1]/ I(\Wb)$ has Krull-dimension $0$. By Lemma \ref{lem: regular sequence Krull dimension}, there exists $\Fb \in \Wb$ such that $\Fb$ is regular in $\ka[\ell_1]$, so $B$ is a quotient of $\ka[\ell_1]/(\Fb)$. Since $\Fb$ is a homogeneous polynomial of degree $d$ in $\ell_1$, it is of the form $c\ell_1^d$ for $c \in \kappa^{\times}$. Then, $\dimk B \leq \dimk \ka[\ell_1]/(\Fb)=d$.
    \item Suppose $h_B(1)=2$ and $h_B(j) \leq j$. Then, by the same reasoning as in $(1)$, we can choose a basis $\ell_3, \dots, \ell_m$ of $(I_B)_1$ and extend them to a basis $(\ell_1, \ell_2, \ell_3, \dots, \ell_m)$ of $R_1$. Let $\Wb$ be the image of the vector space $\langle F_1, \dots, F_m\rangle$ in $R/(\ell_3, \dots, \ell_m) \cong \kappa[\ell_1,\ell_2]$. 
    
    Notice that $\kappa[\ell_1,\ell_2]_{j}$ has dimension $j+1$. Since $h_B(j) \leq j$, there exists $Q \in (I_B)_j$ such that $\bar{Q}$ is non-zero in $\kappa[\ell_1,\ell_2]$, so $\kappa[\ell_1,\ell_2]/(\bar{Q})$ is Cohen-Macaulay of dimension $1$. Consider the quotient $\kappa[\ell_1,\ell_2]/(\Bar{Q},I(\Wb))$ which has Krull dimension $0$. By Lemma \ref{lem: regular sequence Krull dimension}, there exists $\Fb \in \Wb$ such that  $\Fb$ is regular in $\kappa[\ell_1,\ell_2]/(\Bar{Q})$. Let $B':=\kappa[\ell_1,\ell_2]/(\Bar{Q},\Fb)$, so $B'$ has Krull dimension $0$ and $\dimk B' \geq \dimk B$. 

    Now notice that $\dimk B'=jd$. Indeed, we have the exact sequence
    \[0 \to \kappa[\ell_1,\ell_2]/(\Bar{Q})(-d) \xrightarrow{ \cdot \Fb} \kappa[\ell_1,\ell_2]/(\Bar{Q}) \rightarrow B' \to 0\]
Let $C:=\kappa[\ell_1,\ell_2]/(\Bar{Q})$. Then $h_{B'}(n)=h_C(n)-h_C(n-d)$, where $h_C(n)=0$ for $n<0$. Also, \[h_C(n)=\begin{cases}
    j & \text{if } n \geq j \\
    n+1 & \text{if } 0 \leq n \le j-1
\end{cases}\]
So $\dimk B'=\sum_{n=0}^\infty h_{B'}(n)=\sum_{n=0}^\infty (h_C(n)-h_C(n-d))=\sum_{n=N-d+1}^Nh_C(n)=jd$ for sufficiently large $N$. So $\dimk B \leq jd$.
\end{enumerate}
\end{proof}

We now apply Lemma \ref{lem: artinian ring intersection bound} to prove the Propositions \ref{Prop: algebra reg seq 3d-1} and \ref{prop: end bound}.
% \begin{customprop}{\ref{Prop: algebra reg seq 3d-1}}
% Let $A=\kappa[x_1,\ldots,x_m]/I$, where $I$ is a homogeneous
% ideal containing a regular sequence of degrees $(d,\ldots,d)$.
% If $\dim_\kappa A\leq 3d-1$, then $g_A(d)=0$.
% \end{customprop}
\begin{proof}[Proof of Proposition \ref{Prop: algebra reg seq 3d-1}]
Let $A$ be such an Artinian algebra and let $B$ be a graded quotient of $A$. If $h_B(1) \leq 1$, then by Lemma \ref{lem: artinian ring intersection bound}, we have $\dimk B \leq d \leq d+1$. If $h_B(1)=2$ and $h_B(2) \le 2$, then $\dimk B \leq 2d \leq 2d+1$. Hence, by Theorem \ref{thm: 3d-1 artinian}, we get $g_A(d)=0$.
\end{proof}

\begin{proof}[Proof of Proposition \ref{prop: end bound}]
If $d=1$, then $I$ contains all linear forms, so $A=\kappa$ and $g_A(1)=0$. Henceforth assume $d\geq2$.
We separate cases based on $h_A(d)$.

Suppose $h_A(d)\leq d-1$. Then, Corollary~\ref{cor: Gotzman persistence} gives
$h_A(d+i)\leq \max\{d-1-i,0\}$ for every $i\geq 1$.
Hence $g_A(d)\leq \binom{d-1}{2}$.

Suppose $h_A(d)\geq d+1$. Then,  Corollary~\ref{cor: growth condition cor}(1) gives
$h_A(i)\geq i+1$ for $0\leq i\leq d$, so
\[
g_A(d)=\dim_\kappa A-\dim_\kappa A_{\leq d}
\leq d^2+2-\tbinom{d+2}{2}
=\tbinom{d-1}{2}.
\]

It remains to consider $h_A(d)=d$. By Corollary~\ref{cor: growth condition cor}(1),
\[
\dim_\kappa A_{\leq d}
\geq \tsum_{i=0}^{d}\min\{i+1,d\}
=\tbinom{d+2}{2}-1.
\]
If this inequality is strict, then
$\dim_\kappa A_{\leq d}\geq \binom{d+2}{2}$,
and the preceding calculation applies and gives $g_A(d) \leq \binom{d-1}{2}$. Otherwise, equality forces
$h_A(i)=\min\{i+1,d\}$ for every $0\leq i\leq d$.
In particular, $h_A(1)=2$. Applying Lemma~\ref{lem: artinian ring intersection bound} with $B=A$ and $j=d$,
we get $\dim_\kappa A\leq d^2$. Therefore,
\[
g_A(d)
\leq d^2-\big(\tbinom{d+2}{2}-1\big)
=\tbinom{d-1}{2}-1
\leq \tbinom{d-1}{2}. \qedhere
\]
\end{proof}

\section{Bound on $g_A(d)$}\label{Section: Artinian tail}

In this section, we prove Theorem \ref{thm: main inequality}. We state the following lemma.
\begin{lemma}\label{lem: ceil}
Let $d\geq2$.
For a positive integer $n$, let $t_n$ be the smallest integer such that $n\le t_n d+d-t_n+1$. Then $t_n=\ceil{\tfrac{n-(d+1)}{d-1}}$.
\end{lemma}

% \begin{lemma}\label{cor: Gotzman persistence}
% Let $\Gamma\subseteq V(F_1, \dots, F_m) \subseteq \PP^m$. Let $A=\ag$ be the Artinian reduction of $\Gamma$. 

% Then, for $i \geq d$, we have 
% $h_A(i+1)<h_A(i)^{\langle i\rangle}$.

% In particular, if $h_A(i)\le d$, then $h_A(i+1)< h_A(i)$.
% \end{lemma}

% \begin{proof}
% Consider $A'=\kk[x_1, \dots, x_m]/J$, where $J=I(\Gamma)_{\le i}$.
% In particular, $J\supset \langle F_1,\dots, F_m\rangle$.
% so $A'$ is Artinian and $h_{A'}(i)=h_A(i)$ and $h_{A'}(i+1)\ge h_A(i+1)$.
% By Gotzmann persistence applying to $A'$,(cite **), $h_{A'}(i+1)<h_{A'}(i)^{\langle i\rangle}$. So, $h_A(i+1)<h_A(i)^{\langle i\rangle}$. Finally notice that if $i \geq d$ and $h_A(i) \leq d$ then $h_A(i)^{\langle i \rangle}=h_A(i)$.
% \end{proof}
% \begin{theorem}\label{thm: main inequality}
% Let $A=\kappa[x_1,\ldots,x_m]/I$, where $I$ is a homogeneous
% ideal containing a regular sequence of degrees $(d,\ldots,d)$.
% Then, for every $0\leq t\leq d$,
% \begin{equation}\label{eq: main inequality}
% \dim_\kappa A\leq td+d-t+1
% \Longrightarrow
% g_A(d)\leq \tbinom{t-1}{2}.
% \end{equation}
% \end{theorem}
\begin{proof}[Proof of Theorem \ref{thm: main inequality}]
Fix $d\geq2$. We prove the result simultaneously for all
$2\leq t\leq d$ by strong induction on $\dim_\kappa A$.
Proposition~\ref{Prop: algebra reg seq 3d-1} handles
$\dim_\kappa A\leq3d-1$, including the case $t=2$.
Thus, assume $3\leq t\leq d$ and
\[
3d-1<\dim_\kappa A\leq td+d-t+1\leq d^2+1.
\]

For the induction step, we will consider a non-zero linear form $\ell\in A_1$. After fixing a choice of $\ell$, we make the following definitions and observations.
\begin{enumerate}
\item Let $C:= A/\ell A$. Then $C$ is a standard graded Artinian algebra with $\dimk C < \dimk A$.

\item Let $\Ann_A(\ell):=\{a \in A: a\ell=0\}$ be the annihilator of $\ell$ in $A$. Since $\Ann_A(\ell)$ is the kernel of multiplication by $\ell$, we have $A/\Ann_A(\ell) \cong (\ell A)(1)$.
Let $B:=A/\Ann_A(\ell)$. Then $B$ is a standard graded Artinian algebra, and $h_B(i)=\dimk (\ell A)_{i+1}$.

\item Note that we have
\begin{align*}
\dimk A &= \dimk B + \dimk C,
&
g_A(d) &= g_B(d-1) +g_C(d).
\end{align*}

\item Let $s:=h_B(d)$, so $g_B(d-1) = s+g_B(d)$.

\item If $s \geq 2d+1$, then by Corollary~\ref{cor: macaulay lower bound}, we have $\dimk B \geq d^2+s>d^2+1$, which contradicts the fact that $\dimk B \leq \dimk A \leq d^2+1$. Therefore, $s\leq 2d$.

\item If $s=0$, then $h_B(i)=0$ for $i \geq d$, so $g_A(d) = g_C(d)$. Moreover, $\dimk C<\dimk A$, so by induction hypothesis we know that $g_C(d)\leq \binom{t-1}{2}$. Thus we are done in the case when $s=0$, and we can assume $s\geq 1$ for the rest of the proof.

\item Since $1\leq s\leq 2d$, by Corollary~\ref{cor: macaulay lower bound}, we know that
\[
\dimk B - g_B(d-1)
= \dimk B_{\le d-1} 
\geq ds-\tbinom{s}{2}.
\] 

\item Let $t_1$ be the smallest integer such that $\dimk C\leq t_1d+d-t_1+1$. 
If $t_1\geq2$, then $2\leq t_1\leq t\leq d$,
so the induction hypothesis gives
$g_C(d)\leq\binom{t_1-1}{2}$.
If $t_1\leq 2$, then the induction hypothesis gives $g_C(d)=0$.

\item By Lemma \ref{lem: ceil}, we see that
\[
t_1
= \left\lceil\tfrac{\dimk A- \dimk B-(d+1)}{d-1}\right\rceil
\leq \left\lceil\tfrac{(td+d-t+1)- (g_B(d-1) +ds -\binom{s}{2}) -(d+1)}{d-1} \right\rceil
\leq t-s + \left\lceil \tfrac{\binom{s-1}{2} -g_{B}(d-1) }{d-1} \right\rceil.
\]

\item Let $q:=\left\lceil \frac{\binom{s-1}{2} -g_{B}(d-1) }{d-1} \right\rceil$. So $t_1 \leq t-s+q$. Note that $q$ might be positive or negative.

\item 
Since $1\leq s\leq 2d$, we have
\[
q=
\left\lceil\tfrac{\binom{s-1}{2} -g_{B}(d-1) }{d-1} \right\rceil
\leq \left\lceil \tfrac{\binom{s-1}{2}}{d-1} \right\rceil
\leq \left\lceil  \tfrac{(s-1)(2d-2)}{2(d-1)} \right\rceil
= s-1 \leq 2d-1.
\]

\end{enumerate}
We will break the induction step into several cases.

\noindent\textbf{Case 1:} Suppose $h_A(d+1)\leq t-2$. In this case, by Corollary \ref{cor: Gotzman persistence}, we have
\[
g_A(d)
=\dimk A_{\ge d+1}
\le \tsum_{i=1}^{t-2} i
=\tbinom{t-1}{2}.
\]

\noindent\textbf{Case 2:} Suppose $h_A(d+1)\geq t$. We claim that $h_A(d)\ge t+1$. Otherwise, $0<h_A(d)\le t\le d$, so Corollary~\ref{cor: Gotzman persistence} gives $h_A(d+1)<h_A(d)\le t$, a contradiction.

By Corollary~\ref{cor: growth condition cor}(1), we have $h_A(i)\ge\min(i+1,t+1)$ for $0\le i\le d$. Since $t\le d$, it follows that
\[
\dimk A_{\le d}
\ge \tsum_{i=0}^{d}\min(i+1,t+1)
= (d+1)(t+1)-\tbinom{t+1}{2}.
\]
Consequently,
\[
g_A(d)
= \dimk A-\dimk A_{\le d}
\le td+d-t+1-(d+1)(t+1)+\tbinom{t+1}{2}
= \tbinom{t-1}{2}-1.
\]
This proves the desired bound in this case.

Having dealt with Cases 1 and 2, we may therefore assume that $h_A(d+1)=t-1$. In particular, for every nonzero $\ell\in A_1$, we have 
\[
s=\dimk (\ell A)_{d+1}
\leq h_A(d+1)=t-1.
\]

\noindent\textbf{Case 3:} Suppose there exists a linear form $\ell\in A_1$ such that $q\geq 1$.
Then we have
\[
\tbinom{q-1}{2}
=\tfrac{(q-1)(q-2)}{2}
\leq \tfrac{(q-1)(2d-3)}{2}
\leq (q-1)(d-1).
\]
Moreover, by the definition of $q$, we know that $q-1\leq \frac{\binom{s-1}{2}- g_B(d-1) }{d-1}$. Hence, 
\[
\tbinom{q-1}{2}
\leq \tbinom{s-1}{2}-g_B(d-1).
\]
This also implies that $q\leq s$. Also recall that since Case 1 and Case 2 are done, we have the running assumption that $s\leq t-1$.
We therefore see that
\begin{align*}
g_A(d)
&= g_B(d-1) +g_C(d)
\leq g_B(d-1) + \tbinom{t_1-1}{2}
\leq g_B(d-1) + \tbinom{t-s+q-1}{2}\\
&= g_B(d-1) + \tbinom{t-s}{2} +\tbinom{q-1}{2} +(t-s)(q-1)\\
&\leq \tbinom{t-s}{2} +\tbinom{s-1}{2} +(t-s)(q-1)
\leq \tbinom{t-s}{2} +\tbinom{s-1}{2} +(t-s)(s-1)
=\tbinom{t-1}{2}.
\end{align*}

\noindent\textbf{Case 4:} Suppose there exists a linear form $\ell\in A_1$ such that $q\leq 0$ and $t_1\geq2$. Since $t_1\leq t-s+q$, we have $s\leq t-t_1\leq t-2\leq d-2$.

As $B$ is a graded quotient of $A$, Corollary~\ref{cor: Gotzman persistence} applies to $B$. Since $h_B(d)=s\leq d$, its Hilbert function strictly decreases from degree $d$ until it vanishes. Thus $g_B(d-1)\leq s+(s-1)+\cdots+1=\binom{s+1}{2}$.

Recall that the induction hypothesis applied to $C$ gives $g_C(d)\leq \binom{t_1-1}{2}$. Therefore,
\[
g_A(d)
= g_B(d-1)+g_C(d)
\leq \tbinom{s+1}{2}+\tbinom{t_1-1}{2}
\leq \tbinom{t-t_1+1}{2}+\tbinom{t_1-1}{2}
= \tbinom{t-1}{2}-(t_1-2)(t-t_1)
\leq \tbinom{t-1}{2}.
\]
The last inequality follows from $t_1\geq 2$ and $t-t_1\geq s\geq 1$.

\noindent\textbf{Case 5:} Suppose there exists a nonzero linear form $\ell\in A_1$ such that $q\leq 0$ and $\dimk C=2d$. By Proposition~\ref{Prop: algebra reg seq 3d-1}, we have $g_C(d)=0$. 
Recall that after completing Case 1 and Case 2, we have the running assumption that $h_A(d+1)=t-1$. Hence, it follows that $s=h_A(d+1)-h_C(d+1)=t-1$ and $g_A(d)=g_B(d-1)$.

By Corollary~\ref{cor: macaulay lower bound}, we have $\dimk B_{\leq d-1}\geq d(t-1)-\binom{t-1}{2}$. Therefore,
\begin{align*}
g_A(d)
&= g_B(d-1) = \dimk B - \dimk B_{\leq d-1}
= \dimk A-\dimk C-\dimk B_{\leq d-1}\\
&\leq td+d-t+1-2d-d(t-1)+\tbinom{t-1}{2}
= \tbinom{t-1}{2}-t+1
\leq \tbinom{t-1}{2}.
\end{align*}

\noindent\textbf{Case 6:} Suppose there exists a linear form $\ell\in A_1$ such that $q\leq 0$ and $d\leq \dimk C\leq 2d-1$.

We first show that $h_C(d)=0$. Otherwise, Corollary~\ref{cor: Gotzman persistence} implies that $h_C(d-1)\geq 2$. 
Then Corollary~\ref{cor: macaulay lower bound} implies that
\[
\dimk C
\geq h_C(d) +h_C(d-1) +\dimk C_{\leq d-2}
\geq 1+ 2 + 2(d-1) -\tbinom{2}{2}
=2d.
\]
This contradicts the assumption of this case that $\dimk C\leq 2d-1$. Thus, $C_i=0$ for every $i\geq d$. Consequently, $h_B(d-1)=h_A(d)$ and $g_A(d)=g_B(d-1)$.

Since $h_A(d+1)=t-1$ and $t\le d$, Corollary~\ref{cor: Gotzman persistence} implies that $h_A(d)\ge t$. Thus $h_B(d-1)\ge t$, and Corollary~\ref{cor: macaulay lower bound} gives
\[
\dimk B_{\leq d-1} 
= h_B(d-1) +\dimk B_{\leq d-2}
\geq t +t(d-1) -\tbinom{t}{2}
\]
Using $\dimk C\geq d$, we conclude that
\[
g_A(d)
=g_B(d-1)
=\dimk A-\dimk C-\dimk B_{\le d-1}
\leq (td+d-t+1)-d-(dt-\tbinom{t}{2})
=\tbinom{t-1}{2}.
\]

\noindent\textbf{Case 7:} Suppose for every nonzero linear form
$\ell\in A_1$, we have $\dimk C\leq d-1$.

Recall that $h_A(d+1)=t-1$ and $t\geq 3$. Since $A$ is
standard graded, this implies that $h_A(1)\geq 2$.

For every nonzero $\ell\in A_1$, the standard graded algebra
$C=A/\ell A$ satisfies $C_{d-1}=0$. Indeed, otherwise all its
graded components in degrees $0,\ldots,d-1$ would be nonzero,
giving $\dimk C\geq d$. Consequently, multiplication by every
nonzero $\ell\in A_1$ gives a surjective map
$A_i\longrightarrow A_{i+1}$ for every $i\geq d-2$.
In particular, $h_A(i+1)\leq h_A(i)$ in these degrees.

We claim that this inequality is strict whenever $h_A(i)>0$.
Suppose instead that $h_A(i)=h_A(i+1)=r>0$ for some
$i\geq d-2$. Choose two linearly independent forms
$\ell_1,\ell_2\in A_1$, and let $M_1,M_2$ be the matrices of
multiplication by these forms from $A_i$ to $A_{i+1}$, with
respect to fixed bases. Both matrices are invertible.
Thus $\det(M_1+zM_2)$ is a polynomial in $z$ of degree $r$,
with leading coefficient $\det(M_2)\neq 0$. Since $\kappa$
is algebraically closed, this polynomial has a root $z_0\in\kappa$.
Multiplication by the nonzero form $\ell_1+z_0\ell_2$ is
therefore not surjective, which is a contradiction.
We have shown that for $i\geq d-2$, if $h_A(i)>0$, then $h_A(i+1)<h_A(i)$.

Since $h_A(d+1)=t-1$, this strict decrease gives
$h_A(d)\geq t$ and $h_A(d-1)\geq t+1$.
By Corollary~\ref{cor: macaulay lower bound}, it follows that
\[
\dimk A_{\leq d}
=h_A(d)+h_A(d-1) +\dimk A_{\le d-2}
\geq t + (t+1) + (d-1)(t+1) -\tbinom{t+1}{2}.
\]
Consequently,
\[
g_A(d)
=\dimk A-\dimk A_{\leq d}
\leq (td+d-t+1)-\left(t+d(t+1)-\tbinom{t+1}{2}\right)
=\tbinom{t-1}{2}. \qedhere
\]
\end{proof}

\bibliographystyle{plain}
\bibliography{Bibl.bib}
\end{document}